\documentclass[a4paper,12pt]{amsart}
\usepackage{amssymb}
\usepackage{ifthen}
\usepackage{graphicx}
\usepackage{color,xcolor}
\usepackage{array}
\usepackage{caption}
\nonstopmode
\numberwithin{equation}{section}
\newtheorem{thm}{Theorem}
\newtheorem{cor}{Corollary}
\newtheorem{lem}{Lemma}
\newtheorem{prop}{Proposition}

\newtheorem{conj}{Conjecture}
\newtheorem{prob}{Problem}

\theoremstyle{definition}
\newtheorem{defn}{Definition}
\newtheorem{ca}{Case}

\newtheorem{rem}{Remark}

\newenvironment{pf}[1][]{%
 \vskip 1mm
 \noindent
 \ifthenelse{\equal{#1}{}}%
  {{\slshape Proof. }}%
  {{\slshape #1.} }%
 }%
{\qed\medskip}

\newcounter{alphabet}


\newenvironment{Lem}[1][]{\refstepcounter{alphabet}%
\bigskip%
\noindent%
{\bf Lemma \Alph{alphabet}}%
{\bf .} \itshape}{\vskip 8pt}

\newcounter{alphabet2}

\newcommand{\ID}{{\mathbb D}}

\def\be{\begin{equation}}
\def\ee{\end{equation}}

\newcommand{\ben}{\begin{enumerate}}
\newcommand{\een}{\end{enumerate}}

\newcommand{\blem}{\begin{lem}}
\newcommand{\elem}{\end{lem}}
\newcommand{\bthm}{\begin{thm}}
\newcommand{\ethm}{\end{thm}}
\newcommand{\bcor}{\begin{cor}}
\newcommand{\ecor}{\end{cor}}
\newcommand{\beg}{\begin{exam}}
\newcommand{\eeg}{\end{exam}}
\newcommand{\begs}{\begin{examples}}
\newcommand{\eegs}{\end{examples}}
\newcommand{\bdefe}{\begin{defn}}
\newcommand{\edefe}{\end{defn}}
\newcommand{\bprob}{\begin{prob}}
\newcommand{\eprob}{\end{prob}}
\newcommand{\bques}{\begin{ques}}
\newcommand{\eques}{\end{ques}}
\newcommand{\bei}{\begin{itemize}}
\newcommand{\eei}{\end{itemize}}
\newcommand{\bcon}{\begin{conj}}
\newcommand{\econ}{\end{conj}}
\newcommand{\bop}{\begin{op}}
\newcommand{\eop}{\end{op}}

\newcommand{\bas}{\begin{assertion}}
\newcommand{\eas}{\end{assertion}}

\newcommand{\bfa}{\begin{fact}}
\newcommand{\efa}{\end{fact}}

\newcommand{\bca}{\begin{ca}}
\newcommand{\eca}{\end{ca}}

\newcommand{\bst}{\begin{step}}
\newcommand{\est}{\end{step}}

\newcommand{\bsca}{\begin{sca}}
\newcommand{\esca}{\end{sca}}

\newcommand{\bcl}{\begin{cl}}
\newcommand{\ecl}{\end{cl}}

\newcommand{\bmlem}{\begin{mlem}}
\newcommand{\emlem}{\end{mlem}}

\newcommand{\bscl}{\begin{scl}}
\newcommand{\escl}{\end{scl}}

\newcommand{\bcons}{\begin{conjs}}
\newcommand{\econs}{\end{conjs}}

\newcommand{\bprop}{\begin{prop}}
\newcommand{\eprop}{\end{prop}}

\newcommand{\br}{\begin{rem}}
\newcommand{\er}{\end{rem}}
\newcommand{\brs}{\begin{rems}}
\newcommand{\ers}{\end{rems}}
\newcommand{\bo}{\begin{obser}}
\newcommand{\eo}{\end{obser}}
\newcommand{\bos}{\begin{obsers}}
\newcommand{\eos}{\end{obsers}}
\newcommand{\bpf}{\begin{pf}}
\newcommand{\epf}{\end{pf}}
\newcommand{\ba}{\begin{array}}
\newcommand{\ea}{\end{array}}
\newcommand{\beq}{\begin{eqnarray}}
\newcommand{\beqq}{\begin{eqnarray*}}
\newcommand{\eeq}{\end{eqnarray}}
\newcommand{\eeqq}{\end{eqnarray*}}

\newcounter{minutes}
\divide\time by 60
\newcounter{hours}
\multiply\time by 60 \addtocounter{minutes}{-\time}

\begin{document}

\bibliographystyle{amsplain}
\title [On the elliptic harmonic mappings and sense-preserving harmonic mappings] 
{On the elliptic harmonic mappings and sense-preserving harmonic mappings} 

\def\thefootnote{}
\footnotetext{ \texttt{\tiny File:~\jobname .tex,
          printed: \number\day-\number\month-\number\year,
          \thehours.\ifnum\theminutes<10{0}\fi\theminutes}
} \makeatletter\def\thefootnote{\@arabic\c@footnote}\makeatother

\author{Ming-Sheng Liu${}^{~\mathbf{*}}$}
\address{M.-S. Liu, School of Mathematical Sciences, South China Normal University, Guangzhou, Guangdong 510631, China.}
\email{liumsh65@163.com}

\author{Hao XU}
\address{H. XU, School of Mathematical Sciences, South China Normal University, Guangzhou, Guangdong 510631, China. }
\email{136195167@qq.com}

\keywords{harmonic mappings, sense-preserving harmonic mappings, $K$-quasiregular harmonic mappings, $(K,K^{'})$-elliptic harmonic mappings, Landau-Bloch type theorems 
  \\
${}^{\mathbf{*}}$ 
Corresponding author
}


\begin{abstract}
In this paper, we first establish two versions of Landau-Bloch type theorem for $(K,K')$-elliptic harmonic mappings with a bounded minimum distortion. Next, we provide several coefficient estimates and a conjecture for $(K,K')$-elliptic harmonic mappings. Then, we establish three new versions of Landau-Bloch type theorem for sense-preserving harmonic mappings. Finally, we establish two sharp versions of Landau-Bloch type theorem for certain harmonic mappings. These results are sharp in some given cases and improve the related results of different authors.
\end{abstract}

\maketitle

\pagestyle{myheadings}
\markboth{M.-S. Liu and H. Xu}{On the elliptic and sense-preserving harmonic mappings} 

\section{Introduction and Preliminaries}\label{HLP-sec1} 
Suppose that $f(z)=u(z)+i v(z)$ is a twice continuously differentiable function in the unit disk $\mathbb{D}=\{z \in \mathbb{C}:|z|<1\}$. Then $f$ is a harmonic mapping in $\mathbb{D}$ if and only if $f$ satisfies $\Delta f=4 f_{z \bar{z}}=\frac{\partial^2 f}{\partial x^2}+\frac{\partial^2 f}{\partial y^2}=0$ for all
$z=x+i y \in \mathbb{D}$ (cf. \cite{D2004}), where we use the common notations for its formal derivatives:
$$
f_z=\frac{1}{2}\left(f_x-i f_y\right), \quad f_{\bar{z}}=\frac{1}{2}\left(f_x+i f_y\right) .
$$

For such function $f$, we define the maximum and minimum distortions \cite{CGH2000} as follows:
$$
\Lambda_f=\max _{0 \leqslant \theta \leqslant 2 \pi}\left|f_z+e^{-2i \theta} f_{\bar{z}}\right|=\left|f_z\right|+\left|f_{\bar{z}}\right|,
$$
and
$$
\lambda_f=\min _{0 \leqslant \theta \leqslant 2 \pi}\left|f_z+e^{-2i \theta} f_{\bar{z}}\right|=\| f_z|-| f_{\bar{z}}|| .
$$

Since $\mathbb{D}$ is simply connected, $f(z)$ can be written as $f=h+\bar{g}$ with $f(0)=h(0)$, $g(z)$ and $h(z)$ are analytic in $\mathbb{D}$ (cf. \cite{D2004}).
Thus the Jacobian of $f$ is given by
$$
J_f=\left|f_z\right|^2-\left|f_z\right|^2=\left|h^{\prime}\right|^2-\left|g^{\prime}\right|^2.
$$

Notice that $\left|J_f\right|=\lambda_f \cdot \Lambda_f$. It is known \cite{L1936} that a harmonic mapping is locally univalent if and only if  $J_f\neq0$. A harmonic mappings $f$ is said to be sense-preserving if $J_f>0$. Therefore, a sense-preserving harmonic mapping $f$ is said to be $K$-quasiregular harmonic $(K \geq 1)$ on $\mathbb{D}$ if 
$\Lambda_f(z)\leq K \lambda_f(z)$ for all $z\in \mathbb{D}$. There is currently great interest in harmonic mappings due to their relevance for fluid flows (please refer to the discussion in \cite{AC2012,CM2017}).

\textbf{Definition 1.1.} {\rm (\cite{R1993})} Let $\Omega \subset \mathbb{C}$ be a domain. A mapping $f: \Omega \rightarrow \mathbb{C}$ is said to be absolutely continuous on lines, $A C L$ in brief, in a domain $\Omega$ if for every closed rectangle $R \subset \Omega$ with sides parallel to the axes $x$ and $y, f$ is absolutely continuous on almost every horizontal line and almost every vertical line in $R$. Such a mapping $f$ has partial derivatives $f_x$ and $f_y$ a.e. in $\Omega$. Moreover, we say $f \in A C L^2$ if $f \in A C L$ and its partial derivatives are locally $L^2$ integrable in $\Omega$.\\

\textbf{Definition 1.2.} {\rm (\cite{FS1958})} A sense-preserving and continuous mapping $f$ of $\mathbb{D}$ onto $\mathbb{C}$ is said to be $\left(K, K^{\prime}\right)$-elliptic mapping if

(1) $f$ is $A C L^2$ in $\mathbb{D}, J_f \neq 0$ a.e. in $\mathbb{D}$;

(2) $\exists K \geq 1$ and $K^{\prime} \geq 0$ such that $\left\|D_f\right\|^2 \leq K J_f+K^{\prime} \text { a.e. in } \mathbb{D}$.\\

It is easy to see that $\left\|D_f(z)\right\|=\Lambda_f(z)$. Thus, a $\left(K, 0\right)$-elliptic mapping becomes $K$-quasiregular mapping. For more details on elliptic mappings, please refer to (\cite{CPW2020,CPW2021,FS1958,N1953}).

The classical Landau theorem asserts that if $f$ is holomorphic in $\mathbb{D}$ such that $f^{\prime}(0)=1$ and $|f(z)|<M$ for $z\in\mathbb{D}$, then $f$ is univalent in the disk $\mathbb{D}_{\rho_0}=\{z \in \mathbb{C}:|z|<\rho_0\}$, and $f(\mathbb{D}_{\rho_0})$ contains a disk $\mathbb{D}_{R_0}$, where
\begin{eqnarray}
\rho_0=\frac{1}{M+\sqrt{M^2-1}}=M-\sqrt{M^2-1} ~\mbox{ and }~ R_0=M \rho_0^2.
\label{liu11}
\end{eqnarray}
This result is sharp, with the extremal function
\begin{eqnarray}
f_0(z)=M z \left (\frac{1-M z}{M-z}\right ).
\label{liu12}
\end{eqnarray}

The Bloch theorem asserts the existence of a positive constant number b such that if $f$ is an analytic function on the unit disk $\mathbb{D}$ with $f'(0)=1$, then $f(\mathbb{D})$ contains a disk of radius $b$, that is, a disk of radius $b$ which is the univalent image of some region in $\mathbb{D}$. Such a disk is called ``schlicht disk" for $f$. The supremum of all such constants $b$ is called the Bloch constant (see \cite{CGH2000,GK2003}).

In 2000, H.H. Chen {\it et al.} \cite{CGH2000} obtained two versions of Laudau-type theorems for bounded harmonic mappings in the unit disk $\mathbb{D}$.
However, their results are not sharp. Dorff and Nowak \cite{DN2000}, Grigoryan \cite{G2006},
S.L. Chen {\it et al.} \cite{CPR2014,CPW2011C,CPW2011B}, Huang {\it et al.} \cite{H2014, XH2010}, Liu {\it et al.} \cite{L2009C,L2009S,LC2018} and H.H. Chen {\it et al.} \cite{CG2011} improved their results respectively. In particular, the following sharp version of Landau-type theorem for harmonic mappings with bounded maximum distortion was obtained in \cite{H2014,L2009S,LC2018}.

{\bf Theorem A.}\,{\rm (\cite{H2014,L2009S,LC2018})} 
Let $f=h+\overline{g}$ be a harmonic mapping in the unit disk $\mathbb{D}$ such that $f(0)=0$ and $\lambda_f(0)=1$.

\begin{enumerate}
\item If $\varLambda_f(z)\leq 1$ for all $z\in\ID$, then $f$ is univalent in $\mathbb{D}$ and $f(\mathbb{D})$ contains a schlicht disk $\mathbb{D}$,
 and the result is sharp.

\item If $\varLambda_f(z)<\varLambda$ for all $z\in\ID$, then $\varLambda>1$ and $f$ is univalent in the disk $\mathbb{D}_\rho$ with $\rho=\frac{1}{\Lambda}$, and $f(\mathbb{D}_\rho)$ contains a schlicht disk $\mathbb{D}_R$, where
$$
R= \varLambda+(\varLambda^3-\varLambda)\ln\Big (1-\frac{1}{\varLambda^2}\Big).
$$
The result is sharp.
\end{enumerate}

The sharp form of Landau-type theorem for harmonic mappings with $f(0)=0, \lambda_f(0)=1$ and $|f(z)|\leq 1$ in the unit disk was also obtained in \cite{L2009C}. In 2020, Liu {\it et al.} proved the sharp result of Landau-Bloch type theorem for strongly-bounded harmonic mappings when $ M>1 $ in \cite{LLL2020}, and obtained several new sharp versions of Landau-Bloch type theorems of harmonic mappings. We recall one of these results as follows.

{\bf Theorem B.}\, {\rm (\cite[Theorem 3.5]{LLL2020})}\quad Suppose that $M>1$. Let $f(z)$ be a harmonic mapping in the unit disk $\mathbb{D}$ with $f(0)=\lambda_{f}(0)-1=0$, and
$$
f(z)=\sum_{n=1}^{\infty}a_{n}z^{n}+\sum_{n=1}^{\infty}\overline{b_{n}}\overline{z^{n}}
$$
satisfying the following inequality
\begin{eqnarray*}
\sum\limits_{n=2}^{\infty}n(|a_n|+|b_n|)r^{n-1}\leq \frac{(M^2-1)(2Mr-r^2)}{(M-r)^2},\quad 0\leq r\leq \rho_0=M-\sqrt{M^2-1} .
\end{eqnarray*}
Then $f(z)$ is univalent in the disk $\mathbb{D}_{\rho_0}$ and $f(\mathbb{D}_{\rho_0})$ contains a schlicht disk $\mathbb{D}_{R_0}$, where $R_0=M \rho_0^2$. This result is sharp, with $f_0(z)=M z\frac{1-M z}{M-z}$ being an extremal mapping.

Recently, Allu and Kumar \cite{AK2022} established two new versions of Laudau-Bloch type theorem for $(K,K^{'})$-elliptic harmonic mappings and K-quasiregular mappings as follows.

{\bf Theorem C.}\, {\rm (\cite[Theorem 2.5]{AK2022})}  
Let $f=h+\overline{g}$  be a $(K,K^{'})$-elliptic harmonic mapping in $\mathbb{D}$ such that $f(0)=0,\lambda_f(0)=1$ and $\lambda_f(z)\leq \Lambda$ for $z\in\mathbb{D}$. Then $f$ is univalent in the disk $\mathbb{D}_{\rho_1}$ and $f(\mathbb{D}_{\rho_1})$ contains a schlicht disk $\mathbb{D}_{R_1}$,  where
$$
\rho_1=\frac{1}{1+K \Lambda + \sqrt{K'}},\quad
R_1=1+(K\Lambda+\sqrt{K'})\ln\Big(1-\frac{1}{1+K \Lambda + \sqrt{K'}}\Big).
$$

{\bf Theorem D.}\, {\rm (\cite[Theorem 2.7]{AK2022})}
Let $f=h+\overline{g}$  be a K-quasiregular mapping in $\mathbb{D}$ such that $f(0)=0,J_f(0)=1$ and $\lambda_f(z)\leq \Lambda$ for $z\in\mathbb{D}$. Then $f$ is univalent in the disk $\mathbb{D}_{\rho_2}$ and $f(\mathbb{D}_{\rho_2})$ contains a schlicht disk $\mathbb{D}_{R_2}$, where 
$$
\rho_2=\frac{1}{1+K^{3/2} \Lambda },\quad R_2=\frac{1}{\sqrt{K}}+K\Lambda\ln\Big(1-\frac{1}{1+K^{3/2} \Lambda }\Big). 
$$

However, Theorems C and D are not sharp. It is then natural to raise the following problems:

\bprob\label{HLP-prob1}
Can we improve the results in Theorems C and D? Can we establish several sharp results?
\eprob


The paper is organized as follows. In Section \ref{HLP-sec2}, we present statements of our theorems which improve Theorems C and D, and several new versions of Landau-Bloch type theorem and coefficient estimates for sense-preserving harmonic mappings and normalized $(K,K^{'})$-elliptic harmonic mappings, $K$-quasiregular mappings are also provided. In particular, these results are sharp in some given cases. In Section \ref{HLP-sec3}, we state  a couple of lemmas which are needed for the proofs of our main results. In Section \ref{HLP-sec4}, we present the proofs of the main results. In particular, Theorem \ref{HLP-th1}, Corollary \ref{HLP-cr2}, Theorems \ref{HLP-th11}, \ref{HLP-th12}, \ref{HLP-th0} and \ref{HLP-th10} provide an affirmative answer to Problem \ref{HLP-prob1}.

\setcounter{equation}{0}
\section{Statement of Main Results}\label{HLP-sec2}

We first improve Theorems C and D as follows:

\bthm\label{HLP-th1}
Let $f=h+\overline{g}$  be a $(K,K^{'})$-elliptic harmonic mapping defined in the unit disk $\mathbb{D}$ such that $f(0)=0,\lambda_f(0)=1$ and $\lambda_f(z)< \Lambda$ for all $z\in\mathbb{D}$. Then $\Lambda> 1$, $f$ is univalent in the disk $\mathbb{D}_{r_1}$ with
$$r_1=\frac{2}{K\Lambda+\sqrt{K^{2}\Lambda^{2}+4K'}}\geq\frac{1}{K \Lambda + \sqrt{K'}}>\rho_1,$$
and $f(\mathbb{D}_{r_1})$ contains a schlicht disk with radius
$$
\sigma_1=\frac{1}{r_1}+\Big(\frac{1}{r_1^3}-\frac{1}{r_1}\Big)\ln\Big(1-r_1^{2}\Big)>R_1.
$$
When $K=1,\, K'=0$, the radii $r_1=\frac{1}{\Lambda}$ and $\sigma_1=\Lambda+(\Lambda^{3}-\Lambda)\ln\Big(1-\frac{1}{\Lambda^{2}}\Big)$ are sharp.
\ethm

\br 
We remark that Theorem \ref{HLP-th1} improves Theorem C. In order to be more explicit, we compute the approximate values obtained from Theorem \ref{HLP-th1} and the corresponding values obtained from Theorem C for various choices of $K$, $K'$ and $\Lambda$. From Table 1, we can see that, for the same $K$, $K'$ and $\Lambda$, $r_1>\rho_1$, $\sigma_1>R_1$.

\begin{table*}[ht]
\setlength{\belowcaptionskip}{0.2cm}
\centering
\caption{Values of $r_1$, $\rho_1$; $\sigma_1$, $R_1$ in Theorem \ref{HLP-th1}, Theorem C}
\begin{tabular}{ccccccccc}
\hline
  ($K$, $K'$, $\Lambda$) & (1.0,1.1,1.2) & (1.2,1.3,1.4) & (1.4,1.5,1.6) & (1.6,1.7,1.8) & (1.8,1.9,2.0) & (2.0,2.1,2.2)  \\
\hline
  $r_1$ & 0.5530 & 0.4432 & 0.3598 & 0.2956 & 0.2459 & 0.2069 \\
  $\rho_1$ & 0.3078 & 0.2618 & 0.2240 & 0.1929 & 0.1673 & 0.1460  \\
  $\sigma_1$ & 0.3100 & 0.2377 & 0.1882 & 0.1523 & 0.1255& 0.1049  \\
  $R_1$ & 0.1727 & 0.1441  & 0.1214  & 0.1003  & 0.0887  & 0.0768  \\
\hline
\end{tabular}
\end{table*}
\er

Setting $K'=0$ in Theorem \ref{HLP-th1}, we get the following improved version of \cite[Corollary 2.6]{AK2022}. 

\bcor\label{HLP-cr1}
Let $f=h+\overline{g}$  be a $K$-quasiregular mapping defined in the unit disk $\mathbb{D}$ such that $f(0)=0,\lambda_f(0)=1$ and $\lambda_f(z)<\Lambda$ for all $z\in\mathbb{D}$. Then $f$ is univalent in $\mathbb{D}_{r_2}$ 
and $f(\mathbb{D}_{r_2})$ contains a schlicht disk $\mathbb{D}_{\sigma_2}$, where 
$$
r_2=r_2(K,\Lambda)=\frac{1}{K \Lambda}
,\quad \sigma_2=K\Lambda+(K^{3}\Lambda^{3}-K\Lambda)\ln \Big(1-\frac{1}{K^{2}\Lambda^{2}}\Big).
$$
When $K=1$, the result is sharp with an extremal function being $f_1(z)=\Lambda \int_{0}^{z}{\frac{\Lambda z-1}{\Lambda-z}dz}$.
\ecor

\bthm\label{HLP-th3}
Let $f=h+\overline{g}$ be a $(K,K^{'})$-elliptic harmonic mapping in $\mathbb{D}$ such that $f(0)=0,\, J_f(0)=1$ and $\lambda_f(z)< \Lambda$ for all $z\in\mathbb{D}$. Then $\Lambda>  \frac{1}{\sqrt{K+K'}}$, $f$ is univalent in the disk $\mathbb{D}_{r_3}$ and $f(\mathbb{D}_{r_3})$ contains a schlicht disk $\mathbb{D}_{\sigma_3}$, where 
$$
r_3=r_3(K,K',\Lambda)=\frac{2}{\sqrt{K+K'}\, (K\Lambda+\sqrt{K^{2}\Lambda^{2}+4K'})}, 
$$
$$
\sigma_3=\frac{1}{\sqrt{K+K'}}\bigg\{\frac{1}{r_3}+\Big(\frac{1}{r_3^{3}}-\frac{1}{r_3}\Big)\ln\Big(1-r_3^{2}\Big)\bigg\}.
$$
When $K=1,\, K'=0$, the radii $r_3=\frac{1}{\Lambda}$ and $\sigma_3=\Lambda+(\Lambda^{3}-\Lambda)\ln\Big(1-\frac{1}{\Lambda^{2}}\Big)$ are sharp.
\ethm

Setting $K'=0$ in Theorem \ref{HLP-th3}, we get the following result, which improves Theorem D.

\bcor\label{HLP-cr2}
Let $f=h+\overline{g}$  be a K-quasiregular mapping defined on unit disk $\mathbb{D}$ such that $f(0)=0, J_f(0)=1$ and $\lambda_f(z)< \Lambda$ for all $z\in\mathbb{D}$. Then $\Lambda>  \frac{1}{\sqrt{K}}$, $f$ is univalent in the disk $\mathbb{D}_{r_3'}$ and $f(\mathbb{D}_{r_3'})$ contains a schlicht disk $\mathbb{D}_{\sigma_3'}$, where
$$
r_3'=\frac{1}{K^{3/2}\Lambda}>\rho_2,\quad\sigma_3'=K\Lambda+(K^4\Lambda^3-K\Lambda)\ln\Big(1-\frac{1}{K^3\Lambda^2}\Big)>R_2.
$$
When $K=1$, the radii $r_3'=\frac{1}{\Lambda}$ and $\sigma_3'=\Lambda+(\Lambda^{3}-\Lambda)\ln\Big(1-\frac{1}{\Lambda^{2}}\Big)$ are sharp with an extremal function being $f_1(z)=\Lambda \int_{0}^{z}{\frac{\Lambda z-1}{\Lambda-z}dz}$.
\ecor

\br 
From Table 2, we can see that, for the same $K$ and $\Lambda$, $r_3'>\rho_2$, $\sigma_3'>R_2$. 
\begin{table*}[ht]
\setlength{\belowcaptionskip}{0.2cm}
\centering
\caption{Values of $r_3'$, $\rho_3$; $\sigma_3'$, $R_3$ in Corollary \ref{HLP-cr2}, Theorem D}
\begin{tabular}{ccccccccc}
\hline
  ($K$,  $\Lambda$) & (1.0, 1.2) & (1.2, 1.4) & (1.4, 1.6) & (1.6, 1.8) & (1.8, 2.0) & (2.0, 2.2) & (2.2, 2.4) \\
\hline
  $r_3'$ & 0.8333 & 0.5434 & 0.3773& 0.2745 & 0.2070 & 0.1607 & 0.1277 \\
  $\rho_2$ & 0.4545 & 0.3521 & 0.2739 & 0.2154 & 0.1715 & 0.1385& 0.1132 \\
  $\sigma_3'$ & 0.5740 & 0.2768 & 0.1676 & 0.1113 & 0.0783 & 0.0573 & 0.0433 \\
  $R_2$ & 0.2726 &0.1838   &0.1281   &0.0920   &0.0679   &0.0514    &0.0397   \\
\hline
\end{tabular}
\end{table*}
\er

Next, we provide several results and a conjecture of coefficient estimates for normalized $(K,K^{'})$-elliptic harmonic mappings. 

\bthm\label{HLP-th2}
Let $f=h+\overline{g}$ be a $(K,K^{'})$-elliptic harmonic mapping in $\mathbb{D}$ such that $f(0)=0,\, J_f(0)=1$ and $\lambda_f(z)\leq \Lambda$ for all $z\in\mathbb{D}$, and
\begin{eqnarray}
h(z)=\sum\limits_{n=1}^{\infty}{a_nz^n},\quad g(z)=\sum\limits_{n=1}^{\infty}{b_nz^n}.
\label{liu21}
\end{eqnarray}
Then, we have that
\begin{eqnarray*}
\frac{1}{\sqrt{K+K'}}\leq |a_1|+|b_1|\leq \sqrt{K+K'},
\end{eqnarray*}
and
$$
\left|a_n\right|+\left|b_n\right|\leq\frac{(K+K')(K\Lambda+\sqrt{K^{2}\Lambda^{2}+4K'})^2-4}{2n(K+K')(K\Lambda+\sqrt{K^{2}\Lambda^{2}+4K'})},\, n=2,3, \ldots .
$$
When $K=1, K'=0$, the above estimates are sharp for all $n=2,3,\ldots $, with the extremal functions
\begin{eqnarray}
f_n(z)=\Lambda^{2}z-\int_0^z{\frac{\Lambda^{3}-\Lambda}{\Lambda+z^{n-1}}dz}.
\label{liu22}
\end{eqnarray}
\ethm

Applying a similar method as for that used in Theorem \ref{HLP-th2}, we have the following result.

\bcor\label{HLP-cr5}\label{HLP-cr5}
Let $f=h+\overline{g}$ be a $(K,K^{'})$-elliptic harmonic mapping in $\mathbb{D}$ such that $f(0)=0,\, \lambda_f(0)=1$ and $\lambda_f(z)\leq \Lambda$ for all $z\in\mathbb{D}$, and let $h, g$ be given by (\ref{liu21}). Then, we have that
\begin{eqnarray*}
1\leq |a_1|+|b_1|\leq \frac{K+\sqrt{K^2+4K'}}{2},
\end{eqnarray*}
and
$$
\left|a_n\right|+\left|b_n\right|\leq\frac{(K\Lambda+\sqrt{K^{2}\Lambda^{2}+4K'})^2-4}{2n(K\Lambda+\sqrt{K^{2}\Lambda^{2}+4K'})},\, n=2,3, \ldots .
$$
When $K=1, K'=0$, the above estimates are sharp for all $n=2,3,\ldots $, with the extremal functions given by (\ref{liu22}).
\ecor 

Setting $K'=0$ in Theorem \ref{HLP-th2}, we get the following result.

\bcor\label{HLP-cr3}
Let $f=h+\overline{g}$  be a K-quasiregular mapping in $\mathbb{D}$ such that $f(0)=0,J_f(0)=1$ and $\lambda_f(z)\leq \Lambda$ for all $z\in\mathbb{D}$, and let $h, g$ be given by (\ref{liu21}). Then
$$
\frac{1}{\sqrt{K}}\leq |a_1|+|b_1|\leq \sqrt{K},\quad \left|a_n\right|+\left|b_n\right|\leq\frac{K^{3}\Lambda^{2}-1}{n\Lambda K^{2}},\, n=2,3, \ldots .
$$
When $K=1$, the above estimates are sharp for all $n=2,3,\ldots $, with the extremal functions given by (\ref{liu22}).
\ecor

Setting $K'=0$ in Corollary \ref{HLP-cr5}, we get the following result.

\bcor\label{HLP-cr4}
Let $f=h+\overline{g}$  be  a K-quasiregular mapping in $\mathbb{D}$ such that $f(0)=0,\lambda_f(0)=1$ and $\lambda_f(z)\leq \Lambda$ for all $z\in\mathbb{D}$, and let $h, g$ be given by (\ref{liu21}). Then,
$$
1\leq |a_1|+|b_1|\leq K,\quad \left|a_n\right|+\left|b_n\right|\leq\frac{K^2\Lambda^{2}-1}{n\Lambda K}, n=2,3,\ldots .
$$
When $K=1$, the above estimates are sharp for all $n=2,3,\ldots $, with the extremal functions $f_n(z)$ given by (\ref{liu22}).
\ecor

Note that
$$
\frac{K^2\Lambda^{2}-1}{n\Lambda K}=K\frac{\Lambda^{2}-\frac{1}{K^2}}{n\Lambda },
$$
we propose the following conjecture.

\bcon\label{HLP-con2}
Let $F=h+\overline{g}$ be a $K$-quasiregular mapping in unit disk $\mathbb{D}$ such that $F(0)=0,\lambda_F(0)=1$ and $\lambda_F(z)\leq \Lambda$ for all $z\in\mathbb{D}$, and let $h, g$ be given by (\ref{liu21}).  Then,
$$
\left|a_n\right|+\left|b_n\right|\leq K\frac{\Lambda^{2}-1}{n\Lambda },n=2,3,\ldots ,
$$
and the above estimates are sharp for all $n=2,3\ldots$, with the extremal functions:
$$
F_n(z)=\frac{K+1}{2}\bigg(\Lambda^{2}z-\int_{0}^{z}{\frac{\Lambda^{3}-\Lambda}{\Lambda+z^{n-1}}}dz\bigg)+
\overline{\frac{K-1}{2}\bigg(\Lambda^{2}z-\int_{0}^{z}{\frac{\Lambda^{3}-\Lambda}{\Lambda+z^{n-1}}dz}\bigg)}.
$$
\econ


Now, we establish three versions of Landau-Bloch type theorem for sense-preserving harmonic mappings and a new version of Landau-Bloch type theorem for K-quasiregular mapping as follows:

\bthm\label{HLP-th6}
Let $f=h+\overline{g}$  be a sense-preserving harmonic mapping in the unit disk $\mathbb{D}$ such that $f(0)=0,h'(0)-1=g'(0)=0$ and $\lambda_f(z)\leq \lambda$ for all $z\in\mathbb{D}$. Then $\lambda\geq 1$, $f$ is univalent in the disk $\mathbb{D}_{r_4}$ with $r_4=r_4(\lambda)=\frac{1}{2(\sqrt{2}+1)\lambda}$, and $f(\mathbb{D}_{r_4})$ contains a schlicht disk with radius
$$
\sigma_4=\frac{(\sqrt{2}+1)\lambda}{2}\bigg\{1+\Big((\sqrt{2}+1)^2\lambda^2-1\Big)\ln\Big(1-\frac{(\sqrt{2}-1)^2}{\lambda^2}\Big)\bigg\}.
$$
\ethm

\bthm\label{HLP-th11}
 Suppose that $M>1$. Let $f=h+\overline{g}$ be a sense-preserving harmonic mapping in the unit disk $\mathbb{D}$ such that $f(0)=f_z(0)-1=f_{\overline{z}}(0)=0$  and $|h(z)|\leq M$ for all $z\in\mathbb{D}$. Then $f(z)$ is univalent in  $\mathbb{D}_{\rho_0}$ with $\rho_0=\frac{1}{M+\sqrt{M^2-1}}$. The radius $\rho_0$ is sharp with the extremal function $f_0(z)=Mz\frac{1-Mz}{M-z}$.
\ethm

\bthm\label{HLP-th12}
Suppose that $\Lambda>1$. Let $f=h+\overline{g}$ be a sense-preserving harmonic mapping in the unit disk $\mathbb{D}$ such that $f(0)=f_z(0)-1=f_{\overline{z}}(0)=0$ and $|h'(z)|< \Lambda$ for all $z\in\mathbb{D}$. Then $f(z)$ is univalent in $\mathbb{D}_{r_6}$ with $r_6=\frac{1}{\Lambda}$. The radius $r_6$ is sharp with the extremal function $f_1(z)=\int_{0}^{z} {\Lambda\frac{1-\Lambda z}{\Lambda-z}dz}$.
\ethm

\bthm\label{HLP-th7}
Let $f=h+\overline{g}$  be a K-quasiregular mapping defined on unit disk $\mathbb{D}$ such that $f(0)=0,\lambda_f(0)=1$ and $|h(z)|\leq M$ for all $z\in\mathbb{D}$. Then $f$ is univalent in the disk $\mathbb{D}_{r_7}$ and $f(\mathbb{D}_{r_7})$ contains a schlicht disk $\mathbb{D}_{R_7}$, where
$$
r_7=\frac{K+1}{8KM},\quad \sigma_7=\frac{2KM}{K+1}\bigg\{1+\Big(\Big(\frac{4KM}{K+1}\Big)^2-1\Big)\ln\bigg(1-\frac{(K+1)^2}{16K^2M^2}\bigg)\bigg\}.
$$
\ethm


Finally, we establish two sharp versions of Landau-Bloch type theorems for certain harmonic mappings. 

\bthm\label{HLP-th0} Suppose that $M>1$. Let $f=h+\overline{g}$ be a harmonic mapping in the unit disk $\mathbb{D}$ such that $f(0)=0,\, \lambda_f(0)=1$, and let $h(z),g(z)$ satisfy
$$
h(z)=\sum\limits_{n=1}^{\infty}{a_nz^n},\quad g(z)=\sum\limits_{n=1}^{\infty}{b_nz^n}, z\in \mathbb{D},
$$
and for $|z|<\rho_0=\frac{1}{M+\sqrt{M^2-1}}$, we have 
$$
|h'(z)-h'(0)|+|g'(z)-g'(0)|\leq \frac{(M^2-1)(2M|z|-|z|^2)}{(M-|z|)^2}.
$$
Then $f(z)$ is univalent in  $\mathbb{D}_{\rho_0}$ and $f(\mathbb{D}_{\rho_0})$ contains a schlicht disk $\mathbb{D}_{R_0}$, where
$R_0=M \rho_0^2$. This result is sharp, with an extremal function being $f_0(z)=Mz\frac{1-Mz}{M-z}$.
\ethm

\br Theorem \ref{HLP-th0} improves Theorem B, since for $|z|=r<1$, we have
$$
|h'(z)-h'(0)|+|g'(z)-g'(0)|\leq \sum\limits_{n=2}^{\infty}n(|a_n|+|b_n|)r^{n-1}.
$$
\er

\bthm\label{HLP-th10}
Let $f(z)=h(z)+\overline{g(z)}$ be a harmonic mapping in $\mathbb{D}$ with $h(z)$ $g(z)$ given by (\ref{liu21}), such that $f(0)=0, \lambda_f(0)=1$ and $|zf_z(z)+\overline{z}f_{\overline{z}}(z)|<\Lambda$ for all $z\in\mathbb{D}$, 
and $a_nb_n=0$ for $n=2,3,\ldots$. Then $\Lambda\geq 1$ and $f(z)$ 
is univalent in $\mathbb{D}_{r_9}$, and $f(\mathbb{D}_{r_9})$ contains a schlicht disk $\mathbb{D}_{\sigma_9}$, where $r_9=\frac{1}{\Lambda}$, $\sigma_9=\Lambda-\sqrt{\Lambda^2-1}$. The radius $r_9$ is sharp, with the extremal functions $f_1(z)=\Lambda\int_{0}^{z} \frac{1-\Lambda z}{\Lambda-z}dz$. 
Moreover, when $\Lambda=1$, $\sigma_9=1$ is also sharp.
\ethm

\section{Key lemmas}\label{HLP-sec3}
In order to establish our main results, we need the following lemmas.

{\bf Lemma E.}\, {\rm (see, for example, p.35 in \cite{GK2003})}
Let $f(z)=a_0+a_1z+\ldots +a_nz^n+\ldots $ be analytic in $\mathbb{D}$ and let $|f(z)|\leq 1$ for all $z\in\mathbb{D}$. Then $|a_n|\leq1-|a_0|^2$ for $n=1,2,\ldots $.

{\bf Lemma F.}\, {\rm (\cite[Lemma 2.2]{LLL2020})}
Suppose that $\rho>0,\, \sigma>0$, $h(z)$ and $g(z)$ are holomorphic in $\mathbb{D}$ with $h(0)=g(0)=0$. Then, for every complex number $\lambda$ with $\left|\lambda\right|=1$, $f_\lambda(z)=h(z)+\lambda g(z)$ is univalent in $\mathbb{D}_\rho$ and $f_\lambda (\mathbb{D}_\rho)\supseteq \mathbb{D}_\sigma$ if and only if for every complex number $\lambda$ with $\left|\lambda\right|=1$, $F_\lambda(z)=h(z)+\lambda\overline{g(z)}$ is univalent on $\mathbb{D}_\rho$ and $F_\lambda(\mathbb{D}_\rho)\supseteq \mathbb{D}_\sigma$.

{\bf Lemma G.}\, {\rm (see, for example, p.297 in \cite{K2004})}
If $0<x<1$. Then $\ln x>\frac{1}{2}(x-\frac{1}{x})$.

Next, we establish a new lemma, which play a key role in the proof of our main results.

{\bf Lemma H.}\, Suppose that $x\geq 1$. Then $1+x\ln \frac{x}{1+x}<\frac{1}{2(x+0.5)}$.

{\bf Proof:}\quad Let $f(x)=1+x\ln \frac{x}{1+x}-\frac{1}{2x+1}$. Then, a direct computation yields 
$$
f'(x)=\ln\frac{x}{1+x}+\frac{1}{1+x}+\frac{1}{2(x+0.5)^2},
$$
$f'(1)=\frac{1}{2}+\frac{2}{9}-\ln2>0$ and
$$\lim_{x \to +\infty} f'(x)=0.$$
Also, for $x\geq 1$, we have
\begin{eqnarray*}
f''(x)&=&\frac{1}{x(1+x)}-\frac{1}{(1+x)^2}-\frac{1}{(x+0.5)^3}=\frac{1}{x(1+x)^2}-\frac{1}{(x+0.5)^3}\\
&=&-\frac{\frac{1}{2}x^2+\frac{1}{8}(2x-1)}{x(1+x)^2(x+0.5)^3}<0,
\end{eqnarray*}
showing that $f'(x)$ is strictly decreasing on $[1,\infty)$, and thus, we have that $f'(x)>0$ in $[1,\infty)$. This implies that $f(x)$ is strictly increasing in $[1,\infty)$. We note that $f(1)=1-\ln2-\frac{1}{3}<0$ and
$$
\lim_{x \to +\infty} f(x)=0,
$$
and therefore, we obtain that $f(x)<0$ in $[1,\infty)$, and the proof is complete.\hfill $\Box$

\section{Proofs of the main results}\label{HLP-sec4}

\subsection{Proof of Theorem \ref{HLP-th1}}
Since $f=h+\overline{g}$ is a $(K,K^{'})$-elliptic harmonic mapping and $\lambda_f(z)<\Lambda$ for all $z\in\mathbb{D}$, we have that $\Lambda> \lambda_f(0)=1$, and for each $z\in\mathbb{D}$, 
$$
\Lambda_f^2(z)\leq KJ_f(z)+K'\leq K\Lambda_f(z) \lambda_f(z)+K'< K\Lambda\Lambda_f(z) +K',
$$
that is,
$$
 \Lambda_f^2(z)-K\Lambda\Lambda_f(z)-K'< 0.
$$
Hence, for $z\in\mathbb{D}$, we have that
\begin{eqnarray}
\Lambda_f(z)<\frac{K\Lambda+\sqrt{K^{2}\Lambda^{2}+4K'}}{2}. 
\label{liu41}
\end{eqnarray}

Note that $f=h+\overline{g}$ is a harmonic mapping in $\mathbb{D}$ with $f(0)=0,\lambda_f(0)=1$, it follows from Theorem A that $f(z)$ is univalent in the disk $\mathbb{D}_{r_1}$ with
$$r_1=\frac{2}{K\Lambda+\sqrt{K^{2}\Lambda^{2}+4K'}}\geq\frac{1}{K \Lambda + \sqrt{K'}}>\rho_1,$$
and $f(\mathbb{D}_{r_1})$ contains a schlicht disk with radius
$$
\sigma_1=\frac{1}{r_1}+\Big(\frac{1}{r_1^3}-\frac{1}{r_1}\Big)\ln\Big(1-r_1^{2}\Big).
$$

When $K=1,\, K'=0$, we have that $\Lambda(z)=\lambda(z)<\Lambda$ for $z\in\mathbb{D}$, it follows from Theorem A that the radii $r_1=\frac{1}{\Lambda}$ and $\sigma_1=\Lambda+(\Lambda^{3}-\Lambda)\ln\Big(1-\frac{1}{\Lambda^{2}}\Big)$ are sharp.

Finally, we verify that $\sigma_1>R_1$. In fact, it follows from Lemmas G and H that
\begin{eqnarray*}
\sigma_1&\geq& K\Lambda+\sqrt{K'}+\Big((K\Lambda+\sqrt{K'})^{3}-K\Lambda-\sqrt{K'}\Big)\ln\Big(1-\frac{1}{(K\Lambda+\sqrt{K'})^{2}}\Big)\\
&>&K\Lambda+\sqrt{K'}+\Big((K\Lambda+\sqrt{K'})^{3}-K\Lambda-\sqrt{K'}\Big)\cdot\\
&&\cdot\frac{1}{2}\Big(\frac{(K\Lambda+\sqrt{K'})^{2}-1}{(K\Lambda+\sqrt{K'})^{2}}-\frac{(K\Lambda+\sqrt{K'})^{2}}{(K\Lambda+\sqrt{K'})^{2}-1}\Big)
=\frac{1}{2(K\Lambda+\sqrt{K'})}\\
&>&\frac{1}{2(K\Lambda+\sqrt{K'}+0.5)}\\
&>&1+(K\Lambda+\sqrt{K'})\ln\Big(1-\frac{1}{1+K \Lambda + \sqrt{K'}}\Big)=R_1.\hspace{5cm}\Box
\end{eqnarray*}

\subsection{Proof of Theorem \ref{HLP-th3}}
Since $f=h+\overline{g}$ is a $(K,K^{'})$-elliptic harmonic mapping, $J_f(0)=1$ and $\lambda_f(z)<\Lambda$ for all $z\in\mathbb{D}$, we have that $\Lambda> \lambda_f(0)$, and 
$$
\Lambda_f^2(0)\leq KJ_f(0)+K'=K+K', 
$$
that is,
$$
 \Lambda_f(0)\leq \sqrt{K+K'}. 
$$
Hence,
\begin{eqnarray}
\Lambda>\lambda_f(0)=\frac{J_f(0)}{\Lambda_f(0)}\geq\frac{1}{\sqrt{K+K'}}.
\label{liu42}
\end{eqnarray}

Let $F(z)=\frac{f(z)}{\lambda_f(0)}$, then $F(z)$ is a harmonic mapping in $\mathbb{D}$, $\lambda_F(0)=1$, and by (\ref{liu41}), we have
$$\Lambda_F(z)=\frac{\Lambda_f(z)}{\lambda_f(0)}\leq \sqrt{K+K'}\cdot\frac{K\Lambda+\sqrt{K^{2}\Lambda^{2}+4K'}}{2} 
$$
for $z\in \mathbb{D}$. 
So, it follows from Theorem A that $F(z)$ is univalent in the disk $\mathbb{D}_{r_3}$ with
$$
r_3=\frac{2}{\sqrt{K+K'}\, (K\Lambda+\sqrt{K^{2}\Lambda^{2}+4K'})}, 
$$
and $F(\mathbb{D}_{r_3})$ contains a schlicht disk with radius
$$
\sigma_3''=\frac{1}{r_3}+\Big(\frac{1}{r_3^{3}}-\frac{1}{r_3}\Big)\ln\Big(1-r^2_3\Big).
$$

Note that
\begin{eqnarray*}
\lambda_f(0)\sigma_3''\geq  \frac{1}{\sqrt{K+K'}}\bigg\{\frac{1}{r_3}+\Big(\frac{1}{r_3^{3}}-\frac{1}{r_3}\Big)\ln\Big(1-r_3^{2}\Big)\bigg\}=\sigma_3.
\end{eqnarray*}
Hence, $f(z)$ is univalent in the disk $\mathbb{D}_{r_3}$ and $f(\mathbb{D}_{r_3})$ contains a schlicht disk $\mathbb{D}_{\sigma_3}$.

When $K=1,\, K'=0$, we have that $\Lambda(z)=\lambda(z)<\Lambda$ for $z\in\mathbb{D}$, it follows from Theorem A that the radii $r_3=\frac{1}{\Lambda}$ and $\sigma_3=\Lambda+(\Lambda^{3}-\Lambda)\ln\Big(1-\frac{1}{\Lambda^{2}}\Big)$ are sharp. \hfill $\Box$

\subsection{Proof of Corollary \ref{HLP-cr2}}
By Theorem \ref{HLP-th3}, we only need to verify that $\sigma_3'>R_3$. In fact, it follows from Lemmas G and H that
\begin{eqnarray*}
\hspace{1.8cm}\sigma_3'&>&K\Lambda+(K^4\Lambda^3-K\Lambda)\cdot\frac{1}{2}\Big(\frac{K^{3}\Lambda^{2}-1}{K^{3}\Lambda^{2}}-\frac{K^{3}\Lambda^{2}}{K^{3}\Lambda^{2}-1}\Big)=\frac{1}{2K^2\Lambda}\\
&>&\frac{1}{\sqrt{K}}\cdot\frac{1}{2(K^{3/2}\Lambda+0.5)}\\
&>&\frac{1}{\sqrt{K}}\bigg(1+K^{3/2}\Lambda\ln\Big(1-\frac{1}{1+K^{3/2}\Lambda}\Big)\bigg)=R_3.\hspace{3.6cm}\Box
\end{eqnarray*}

\subsection{Proof of Theorem \ref{HLP-th2}}
Since $f=h+\overline{g}$ is a $(K,K^{'})$-elliptic harmonic mapping, $J_f(0)=1$ and $\lambda_f(z)<\Lambda$ for all $z\in\mathbb{D}$, it follows from (\ref{liu41}) and (\ref{liu42}) that
\begin{eqnarray*}
||a_1|-|b_1||&=&\lambda_f(0)\geq\frac{1}{\sqrt{K+K'}}, \\
\Lambda_f(z)&<&\Lambda_1:=\frac{K\Lambda+\sqrt{K^{2}\Lambda^{2}+4K'}}{2},\, z\in\mathbb{D}.
\end{eqnarray*}
Hence
\begin{eqnarray*}
\frac{1}{\sqrt{K+K'}}\leq\lambda_f(0)\leq |a_1|+|b_1|=\frac{J_f(0)}{\lambda_f(0)}=\frac{1}{\lambda_f(0)}\leq \sqrt{K+K'}.
\end{eqnarray*}

Let $H(z)=h(z)+e^{i\alpha}g(z),\, 0\leq\alpha\leq2\pi$, then $H$ is holomorphic in $\mathbb{D}$ and for $z\in \mathbb{D}$,
$$
|H'(z)|=\bigg|\sum\limits_{n=1}^{\infty}{n(a_n+e^{i\alpha}b_n)z^{n-1}}\bigg|\leq \Lambda_f(z)<\Lambda_1. 
$$

By Lemma E, we have that
$$
\Big|\frac{n(a_n+e^{i\alpha}b_n)}{\Lambda_1}\Big|\leq 1-\Big|\frac{a_1+e^{i\alpha}b_1}{\Lambda_1}\Big|^2\leq1-\Big|\frac{|a_1|-|b_1|}{\Lambda_1}\Big|^2
\leq 1-\Big|\frac{1}{\sqrt{K+K'}\Lambda_1}\Big|^2.
$$
Then, it follows from the arbitrary of $\alpha$ that
$$
|a_n|+|b_n|\leq\frac{(K+K')\Lambda_1^2-1}{n(K+K')\Lambda_1}=\frac{(K+K')(K\Lambda+\sqrt{K^{2}\Lambda^{2}+4K'})^2-4}{2n(K+K')(K\Lambda+\sqrt{K^{2}\Lambda^{2}+4K'})},\quad n=2,3,\ldots .
$$

When $K=1,\, K'=0$, the above estimates are sharp for all $n=2,3,\ldots$, with the extremal functions
$$f_n(z)=\Lambda^{2}z-\int_0^z{\frac{\Lambda^{3}-\Lambda}{\Lambda+z^{n-1}}dz}= z+\frac{\Lambda^2-1}{n\Lambda}z^n+\sum\limits_{k=2}^{\infty}
\frac{(-1)^{k+1}(\Lambda^2-1)}{(kn-k+1)\Lambda ^k} z^{kn-k+1}.\eqno\Box$$

\subsection{Proof of Theorem \ref{HLP-th6}}
Since $f=h+\overline{g}$ is a sense-preserving harmonic mapping defined on unit disk $\mathbb{D}$ such that $f(0)=0, h'(0)-1=g'(0)=0$, we have that $|h'(z)|>|g'(z)|$ for $z\in\mathbb{D}$, this implies that $w(z)=\frac{g'(z)}{h'(z)}$ is holomorphic in $\mathbb{D}$, $ w(0)=0$ and $|w(z)|<1$ for $z\in\mathbb{D}$.

By Schwarz's lemma, we have that $|w(z)|\leq |z|$, i.e., $|g'(z)|\leq|z||h'(z)|$.

Since $\lambda_f(z)=|h'(z)|-|g'(z)|\leq\lambda$, we have that
$$
|h'(z)|\leq\frac{\lambda}{1-|z|}\mbox{ and } |g'(z)|\leq\frac{\lambda |z|}{1-|z|}\mbox{  for  }z\in\mathbb{D}.
$$
So,
$$
\Lambda_f(z)=|h'(z)|+|g'(z)|\leq\frac{\lambda (1+|z|)}{1-|z|}.
$$

Let $P(z)=(\sqrt{2}+1)f(\frac{z}{\sqrt{2}+1}),\, z\in\mathbb{D}$, a simple computation yields that $\lambda_P(0)=\lambda_f(0)=1$, and
$$
\Lambda_P(z)=\Lambda_f(\frac{z}{(\sqrt{2}+1)})\leq\frac{\lambda (1+|z/(\sqrt{2}+1)|)}{1-|z/(\sqrt{2}+1)|}\leq (\sqrt{2}+1)\lambda . 
$$

Therefore, by Theorem A, we obtain that $P(z)$ is univalent in the disk $\mathbb{D}_{r_4'}$ and $F(\mathbb{D}_{r_4'})$ contains a schlicht disk $\mathbb{D}_{\sigma_4'}$, where
$$
r_4'=\frac{1}{(\sqrt{2}+1)\lambda},\, \sigma_4'=(\sqrt{2}+1)\lambda+\Big((\sqrt{2}+1)^{3}\lambda^3-(\sqrt{2}+1)\lambda\Big)\ln\Big(1-\frac{1}{(\sqrt{2}+1)^2\lambda^2}\Big).
$$
Hence, $f$ is univalent in $\mathbb{D}_{r_4}$ and $f(\mathbb{D}_{r_4})$ contains a schlicht disk $\mathbb{D}_{\sigma_4}$, where
\begin{eqnarray*}
r_4&=&\frac{r_4'}{2}=\frac{1}{2(\sqrt{2}+1)\lambda},\quad \\ \sigma_4&=&\frac{\sigma_4'}{2}=\frac{(\sqrt{2}+1)\lambda}{2}\bigg\{1+\Big((\sqrt{2}+1)^2\lambda^2-1\Big)\ln\Big(1-\frac{(\sqrt{2}-1)^2}{\lambda^2}\Big)\bigg\}.
\hspace{2cm}\Box
\end{eqnarray*}

\subsection{Proof of Theorem \ref{HLP-th11}}
We adopt a similar method as for that used in \cite[Theorem 2]{CG2011}. Let $F(z)=\int_{0}^{z}{f_z(z)dz}$ for $z\in\mathbb{D}$, then it follows from $f_z(z)$ is holomorphic in $\mathbb{D}$ that $F(z)$ is holomorphic in $\mathbb{D}$. Since $f(0)=f_z(0)-1=f_{\overline{z}}(0)=0$ and $|h(z)|\leq M$ for $z\in\mathbb{D}$, we have that
$$F(0)=0, \quad F'(0)=f_z(0)=1, \quad |F(z)|=\bigg|\int_{0}^{z}{f_z(z)dz}\bigg|=\bigg|\int_{0}^{z}{h'(z)dz}\bigg|=|h(z)|\leq M.$$

By the classical Landau theorem, we have that $F(z)$ is univalent in $\mathbb{D}_{\rho_0}$ 
and $F(\mathbb{D}_{\rho_0})$ contains a disk $\mathbb{D}_{R_0}$, where $\rho_0$ is given by (\ref{liu11}) and $R_0=M \rho_0^2$.

Since $f(z)$ is a sense-preserving harmonic mapping in $\mathbb{D}$, we have that $|f_z(z)|>|f_{\overline{z}}(z)|$ for all $z\in\mathbb{D}$, this implies that $f_z(z)\neq 0$ for all $z\in\mathbb{D}$. 
Thus, the function $\frac{\overline{f_{\overline{z}}}}{f_z}$ is holomorphic and has its modulus less than 1 on $\mathbb{D}$. 

Let $\Delta=F(\mathbb{D}_{\rho_0})$, then $F(z)$ is biholomorphic from $\mathbb{D}_{\rho_0}$ onto $\Delta$. By applying Schwarz's lemma to $F(z)$, it is easy to verify that $\Delta\subset \mathbb{D}$, therefore, the composition mapping $g(\zeta)=f\circ F^{-1}(\zeta)$ is a harmonic mapping in $\Delta\, (\subset \mathbb{D})$, and
$$g_\zeta=\frac{f_z(z)}{F'(z)}=1,\quad |g_{\overline{\zeta}}|=\bigg|\frac{f_{\overline{z}}}{\overline{f_z}}\bigg|<1 $$
for $\zeta=F(z) \in \Delta $. Thus, for any $\zeta_1 ,\zeta_2 \in \Delta$ with $\zeta_1\neq\zeta_2$, we have
$$
|g(\zeta_1)-g(\zeta_2)|\geq \bigg|\int_{\overline{\zeta_1 ,\zeta_2}}{g_\zeta d\zeta}\bigg|-\bigg|\int_{\overline{\zeta_1 ,\zeta_2}}{g_{\overline{\zeta}} d\overline{\zeta}}\bigg|>0.
$$
This shows that $g$ is univalent on $\Delta$ and, consequently,  $f(z)$ is univalent on $\mathbb{D}_{\rho_0}$.

Finally, we prove the sharpness of $\rho_0$. To this end, we consider the harmonic mapping
$$f(z)=h(z)=Mz\frac{1-Mz}{M-z}.$$

Note that $f(0)=0, f_{\overline{z}}(z)\equiv 0$ and for $z\in \mathbb{D}$,
$$ |h(z)|=|f(z)|=\Big|Mz\frac{1-Mz}{M-z}\Big|< M, \quad f'(z)=M^2\frac{(z-M-\sqrt{M^2-1})(z-\rho_0)}{(M-z)^2}, $$
we obtain that $f_z(0)-1=f_{\overline{z}}(0)=0$, $|f_z(z)|>|f_{\overline{z}}(z)|\equiv 0$ for $z\in \mathbb{D}\backslash\{\rho_0\}$, and $f'(\rho_0)=0$, which proves the sharpness of $\rho_0$. \hfill $\Box$

\subsection{Proof of Theorem \ref{HLP-th12}}
As in the proof of Theorem \ref{HLP-th11}, the function $F(z):=\int_{0}^{z}{f_z(z)dz}$ is holomorphic in $\mathbb{D}$.
Since $f(0)=f_z(0)-1=f_{\overline{z}}(0)=0$ and $|h'(z)|<\Lambda$ for $z\in\mathbb{D}$, we have that
$$F(0)=0, \quad F'(0)=f_z(0)=1, \quad |F'(z)|=|f_z(z)|=|h'(z)|<\Lambda\quad\mbox{ for } z\in\mathbb{D}.$$

By Theorem A, we obtain that $F(z)$ is univalent in the disk $\mathbb{D}_{r_6}$ with $r_6=\frac{1}{\Lambda}$.

Note that
$$
|F(z)|=\Big|\int_{0}^{z}h'(t)dt\Big|\leq \int_{0}^{z}|h'(t)||dt|\leq \Lambda |z|<\Lambda\quad\mbox{ for } z\in\mathbb{D}.
$$

By using an analogous proof of Theorem \ref{HLP-th11} and the extremal function $f_1(z)=h(z)=\int_{0}^{z}{\Lambda\frac{1-\Lambda z}{\Lambda-z}dz}$, we may prove that $f(z)$ is univalent on $\mathbb{D}_{r_6}$, and the radius $r_6$ is sharp.\hfill $\Box$




\subsection{Proof of Theorem \ref{HLP-th7}}
Since $f$ is a $K$-quasiregular mapping in $\mathbb{D}$, we have that
$$
|f_z|+|f_{\overline{z}}|\leq K(|f_z|-|f_{\overline{z}}|),
$$
that is,
$$
|f_{\overline{z}}|\leq\frac{K-1}{K+1}|f_z|.
$$

Because $h(0)=0$ and $|h(z)|<M$ for $z\in\mathbb{D}$, by Cauchy's estimates theorem, we have
$$
|f_z(z)|=|h'(z)|\leq\frac{M}{ 1-|z|}
$$
for $z\in\mathbb{D}$. Then
$$
\Lambda_f(z)=|f_z|+|f_{\overline{z}}|\leq\frac{2KM}{ (K+1)(1-|z|)}.
$$

Let $P(z)=2f(\frac{z}{2}),\, z\in\mathbb{D}$. Then $\lambda_P(0)=\lambda_f(0)=1$, and
$$
\Lambda_P(z)=\Lambda_f(\frac{z}{2})\leq\frac{2KM}{ (K+1)(1-|z/2|)}\leq\frac{4KM}{K+1}.
$$

Therefore, by Theorem A, we obtain that $P(z)$ is univalent in the disk $\mathbb{D}_{r_7'}$ and $F(\mathbb{D}_{r_7'})$ contains a schlicht disk $\mathbb{D}_{\sigma_7'}$, where
$$
r_7'=\frac{K+1}{4KM},\quad \sigma_7'=\frac{4KM}{K+1}+\Big(\Big(\frac{4KM}{K+1}\Big)^3-\frac{4KM}{K+1}\Big)\ln\bigg(1-\frac{(K+1)^2}{16K^2M^2}\bigg).
$$
Hence, $f$ is univalent in $\mathbb{D}_{r_7}$ and $f(\mathbb{D}_{r_7})$ contains a schlicht disk $\mathbb{D}_{\sigma_7}$, where
\begin{eqnarray*}
r_7&=&\frac{r_7'}{2}=\frac{K+1}{8KM},\quad \\ \sigma_7&=&\frac{\sigma_7'}{2}=\frac{2KM}{K+1}\bigg\{1+\Big(\Big(\frac{4KM}{K+1}\Big)^2-1\Big)\ln\bigg(1-\frac{(K+1)^2}{16K^2M^2}\bigg)\bigg\}.\hspace{2cm}\Box
\end{eqnarray*}

\subsection{Proof of Theorem \ref{HLP-th0}}
Since $\lambda_f(0)=1$, we have that $||a_1|-|b_1||=\lambda_f(0)=1$.

For every complex number $\lambda$ with $|\lambda|=1$, let
$$f_\lambda(z)=h(z)+\lambda g(z)=\sum\limits_{n=1}^{\infty}{(a_n+\lambda b_n)z^n}.$$

Since $|h'(z)-h'(0)|+|g'(z)-g'(0)|\leq \frac{(M^2-1)(2M|z|-|z|^2)}{(M-|z|)^2} $ for $|z|<\rho_0$, we have
\begin{eqnarray*}
\Big|\sum\limits_{n=2}^{\infty}{n(a_n+\lambda b_n)z^{n-1}}\Big|&\leq& |h'(z)-h'(0)|+|g'(z)-g'(0)|\\ 
& \leq & \frac{(M^2-1)(2M|z|-|z|^2)}{(M-|z|)^2}\medskip\\
& < & 
\frac{(M^2-1)(2M\rho_0-\rho_0^2)}{(M-\rho_0)^2}
\end{eqnarray*} for $|z|<\rho_0$.

We first prove $f_\lambda(z)$ is univalent in $\mathbb {D}_{\rho_0}$. To this end, for any $z_1\neq z_2\in\mathbb{D}_{\rho_0}$, we have
\begin{eqnarray*}
|f_\lambda(z_1)-f_\lambda(z_2)|& = & \bigg|\int_{[z_1,z_2]}{\sum\limits_{n=1}^{\infty}{n(a_n+\lambda b_n)z^{n-1}}}dz\bigg|\medskip\\
& \geq & \bigg|\int_{[z_1,z_2]}{(a_1+\lambda b_1)dz}\bigg|-\bigg|\int_{[z_1,z_2]}{\sum\limits_{n=2}^{\infty}{n(a_n+\lambda b_n)z^{n-1}}dz}\bigg|\medskip\\
& \geq & |z_1-z_2|||a_1|-|b_1||-\int_{[z_1,z_2]}{\bigg|\sum\limits_{n=2}^{\infty}{n(a_n+\lambda b_n)z^{n-1}}\bigg||dz|}
\medskip\\
& > & |z_1-z_2|-\int_{[z_1,z_2]}{\frac{(M^2-1)(2M\rho_0-\rho_0^2)}{(M-\rho_0)^2}|dz|}\medskip\\
& \geq & |z_1-z_2|\bigg(1-\frac{(M^2-1)(2M\rho_0-\rho_0^2)}{(M-\rho_0)^2}\bigg)=0.
\end{eqnarray*}
Thus, we have that $f_\lambda(z_1)\neq f_\lambda(z_2)$, which proves the univalence of $f_\lambda(z)$ in $\mathbb {D}_{\rho_0}$.

Next, noticing that $f(0)=0$, for any $z'\in\partial\mathbb{D}_{\rho_0}$, we have
\begin{eqnarray*}
|f_\lambda(z')|& = & |\Sigma_{n=1}^\infty{(a_n+\lambda b_n)z'^n}| \geq \lambda_f(0) \rho_0 - |\Sigma_{n=2}^\infty{(a_n+\lambda b_n)z'^n}|\medskip\\
& \geq & \rho_0-\int_{[0,z']}{|\Sigma_{n=2}^\infty{n(a_n+\lambda b_n)z^{n-1}|}}|dz|\medskip\\
& \geq & \rho_0-\int_{[0,z']}{\frac{(M^2-1)(2M|z|-|z|^2)}{(M-|z|)^2}|dz|}\medskip\\
& =& \rho_0-\int_{0}^{1}{\frac{(M^2-1)(2Mt\rho_0-t^2\rho_0^2)}{(M-t\rho_0)^2}d(t \rho_0)}\medskip\\
& =&\rho_0-\frac{(M^2-1)\rho_0^2}{M-\rho_0}=M\rho_0\frac{1-M\rho_0}{M-\rho_0}=M\rho_0^2=R_0.
\end{eqnarray*}
Hence, $f_\lambda(\mathbb{D}_{\rho_0})$ contains a schlicht disk $\mathbb{D}_{R_0}$.
By Lemma F, we obtain that $f(z)$ is univalent on $\mathbb{D}_{\rho_0}$ and $f(\mathbb{D}_{\rho_0})$ contains a schlicht disk $\mathbb{D}_{R_0}$.

Finally, we prove the sharpness of $\rho_0$ and $R_0$. We consider the harmonic mapping $f_0(z)=Mz\frac{1-Mz}{M-z}$. Note that $f_0(0)=0$,
$\lambda_{f_0}(0)=1$ and
$$
f_0(z)=Mz\frac{1-Mz}{M-z}=z-\sum\limits_{n=2}^{\infty}{\frac{M^2-1}{M^{n-1}}z^n},z\in \mathbb{D},
$$
direct computation yields
$$
\bigg|\sum\limits_{n=2}^{\infty}{n\frac{M^2-1}{M^{n-1}}z^{n-1}}\bigg|\leq \sum\limits_{n=2}^{\infty}{n\frac{M^2-1}{M^{n-1}}|z|^{n-1}}=\frac{(M^2-1)(2M|z|-|z|^2)}{(M-|z|)^2},\quad 0\leq r<1,
$$
thus it follows from the above results that $f_0(z)$ is univalent in $\mathbb{D}_{\rho_0}$ and $f(\mathbb{D}_{\rho_0})$ contains a schlicht disk $\mathbb{D}_{R_0}$.

Since
$$f'_0(z)=M^2\frac{(z-M-\sqrt{M^2-1})(z-\rho_0)}{(M-z)^2},$$
we obtain that $f'_0(\rho_0)=0$, which proves the sharpness of $\rho_0$.

Noticing that $f_0(0)=0$, $z'=\rho_0\in\partial\mathbb{D}_{\rho_0}$, we have
$$|f_0(z')-f_0(0)|=|f_0(z')|=M\rho_0\frac{1-M\rho_0}{M-\rho_0}=M\rho_0^2=R_0.$$
Hence, the radius $R_0$ is also sharp. This completes the proof. \hfill $\Box$

\subsection{Proof of Theorem \ref{HLP-th10}}
Let $z=re^{i\theta},\, r\in [0,1)$, $\theta \in [0,2\pi)$, then
$$
f(re^{i\theta})=\sum\limits_{n=1}^{\infty}{a_nr^ne^{in\theta}}+\sum\limits_{n=1}^{\infty}{\overline{b_n}r^ne^{-in\theta}}.
$$

Since $|zf_z(z)+\overline{z}f_{\overline{z}}(z)|<\Lambda$ for all $z\in\mathbb{D}$, 
by Parseval's identity, we have that
$$
\sum\limits_{n=1}^{\infty}{n^2(|a_n|^2+|b_n|^2)r^{2n}}=\frac{1}{2\pi}\int_0^{2\pi}{|re^{i\theta}f_z(re^{i\theta})+\overline{re^{i\theta}}f_{\overline{z}}(re^{i\theta})|^2 d\theta}\leq \Lambda^2.
$$

Notice that $\lambda_f(0)=||a_1|-|b_1||=1$, 
let $r\rightarrow 1^-$, we have that $\Lambda\geq 1$, and 
$$\sum\limits_{n=2}^{\infty}{n^2(|a_n|^2+|b_n|^2)}=\sum\limits_{n=1}^{\infty}{n^2(|a_n|^2+|b_n|^2)}-(|a_1|^2+|b_1|^2)\leq \Lambda^2-1.$$

Since $a_nb_n=0$ for $n=2,3,\ldots$, we have that
\begin{eqnarray*}
\sum\limits_{n=2}^{\infty}{n^2(|a_n|+|b_n|)^2}=\sum\limits_{n=2}^{\infty}{n^2(|a_n|^2+|b_n|^2)}\leq \Lambda^2-1.
\end{eqnarray*}

For any complex number $\lambda$ with $|\lambda|=1$, let
$$
f_\lambda(z)=h(z)+\lambda g(z)=\int_{0}^{z}{\sum\limits_{n=1}^{\infty}{n(a_n+\lambda b_n)z^{n-1}}dz}.
$$

We first prove $f_\lambda(z)$ is univalent in the $\mathbb {D}_{r_9}$ with $r_9=\frac{1}{\Lambda}$. To this end, for any $z_1, z_2\in\mathbb{D}_{r}$ with $z_1\neq z_2$ and $0<r<r_9$, we have
\begin{eqnarray*}
|f_\lambda(z_1)-f_\lambda(z_2)|& = & \bigg|\int_{[z_1,z_2]}{\sum\limits_{n=1}^{\infty}{n(a_n+\lambda b_n)z^{n-1}}}dz\bigg|\medskip\\
& \geq & \bigg|\int_{[z_1,z_2]}{(a_1+\lambda b_1)dz}\bigg|-\bigg|\int_{[z_1,z_2]}{\sum\limits_{n=2}^{\infty}{n(a_n+\lambda b_n)z^{n-1}}dz}\bigg|\medskip\\
& \geq & |z_1-z_2|||a_1|-|b_1||-\int_{[z_1,z_2]}{\bigg|\sum\limits_{n=2}^{\infty}{n(a_n+\lambda b_n)z^{n-1}}\bigg||dz|}\\
& \geq & |z_1-z_2|\bigg(1-\sum\limits_{n=2}^{\infty}{n(|a_n|+|b_n|)r^{n-1}}\bigg)\medskip\\
& \geq & |z_1-z_2|\bigg(1-\sqrt{\sum\limits_{n=2}^{\infty}{n^2(|a_n|+|b_n|)^2}  \sum\limits_{n=2}^{\infty}{r^{2n-2}}}\bigg)\medskip\\
& \geq & |z_1-z_2|\bigg(1-\sqrt{(\Lambda^2-1)  \sum\limits_{n=2}^{\infty}{r^{2n-2}} }\bigg)\medskip\\
& = & |z_1-z_2|\bigg(1-\frac{r\sqrt{\Lambda^2-1}}{\sqrt{1-r^2}}\bigg)> 0,
\end{eqnarray*}
which implies that $f_\lambda(z_1)\neq f_\lambda(z_2)$, this proves the univalence of $f_\lambda(z)$ in the $\mathbb {D}_{r_9}$.

Next, noticing that $f_\lambda(0)=0$, for any $z'\in\partial\mathbb{D}_{r_9}$, we have
\begin{eqnarray*}
|f_\lambda(z^{'})|& = & \bigg|\int_{[0,z']}{\sum\limits_{n=1}^{\infty}{n(a_n+\lambda b_n)z^{n-1}}}dz\bigg|\medskip\\
\end{eqnarray*}\begin{eqnarray*}
& \geq &   |z'|||a_1|-|b_1||-\int_{[0,z']}{\bigg|\sum\limits_{n=2}^{\infty}{n(a_n+\lambda b_n)z^{n-1}}\bigg||dz|}
\medskip\\
& \geq & r_9-\int_{0}^{r_9}{\Sigma_{n=2}^\infty{n(|a_n|+| b_n|)r^{n-1}dr}}
\medskip\\
& \geq & r_9-\int_{0}^{r_9}{ \sqrt{\sum\limits_{n=2}^{\infty}{n^2(|a_n|+|b_n|)^2}  \sum\limits_{n=2}^{\infty}{r^{2n-2}}}\, dr}
\medskip\\
& \geq & r_9-\sqrt{\Lambda^2-1}\int_{0}^{r_9}\frac{r}{\sqrt{1-r^2}}\,  dr\medskip\\
& =&r_9+\sqrt{\Lambda^2-1}(\sqrt{1-r_9^2}-1)=\Lambda-\sqrt{\Lambda^2-1}=\sigma_9.
\end{eqnarray*}
Hence, $f_\lambda(\mathbb{D}_{r_9})$ contains a schlicht disk $\mathbb{D}_{\sigma_9}$. By Lemma F, we have that $f(z)$ is univalent in $\mathbb{D}_{r_9}$ and $f(\mathbb{D}_{r_9})$ contains a schlicht disk $\mathbb{D}_{\sigma_9}$.

Now, we prove the sharpness of $r_9$. We consider the harmonic mapping
$$f_1(z)=\int_{0}^{z}{\Lambda\frac{1-\Lambda z}{\Lambda-z}dz}.$$

It is easy to verify that $f_1(z)$ satisfy all of the hypothesis of Theorem \ref{HLP-th10}.
Then it follows from the above results that the harmonic mapping $f_1(z)$ is univalent in $\mathbb{D}_{r_9}$ and $f(\mathbb{D}_{r_9})$ contains a schlicht disk $\mathbb{D}_{\sigma_9}$. By $f_1'(r_9)=0$, we get that the radius $r_9$ is sharp.

Finally, when $\Lambda=1$, it is obvious that $\sigma_9=r_9=1$ are sharp with the extremal mapping $f_2(z)=z$. \hfill $\Box$

\subsection*{Acknowledgments}
This research is supported by Guangdong Natural Science Foundations (Grant No. 2021A1515010058).

\subsection*{Conflict of Interests}
The authors declare that they have no conflict of interest, regarding the publication of this paper.

\subsection*{Data Availability Statement}
The authors declare that this research is purely theoretical and does not associate with any data.

\end{document}